\documentclass{birkart} 
\usepackage{amssymb,amsmath,exscale,enumerate}

\usepackage{amsthm}

\usepackage{amsopn}
\usepackage{array}
\usepackage{verbatim}

\DeclareMathOperator{\supp}{supp}

\DeclareMathOperator{\Id}{Id}
\DeclareMathOperator{\Tr}{Tr}
     
     \newcommand{\EE}{\mathbb{E}}
     
     \newcommand{\NN}{\mathbb{N}}
     \newcommand{\PP}{\mathbb{P}}
     
     \newcommand{\RR}{\mathbb{R}}
     \newcommand{\ZZ}{\mathbb{Z}}
\newcommand{\cD}{\mathcal{D}}
\newcommand{\cP}{\mathcal{P}}
\newcommand{\cU}{\mathcal{U}}
\newcommand{\co}{\mathfrak{o}}
\newcommand{\be}{\begin{equation}}
\newcommand{\ee}{\end{equation}}
\newcommand{\la}{\langle}
\newcommand{\ra}{\rangle}
\newcommand{\D}[1]{\Delta^{\scriptscriptstyle #1}}
\newcommand{\N}[1]{N^{\scriptscriptstyle #1}}
\newcommand{\Si}[1]{\Sigma^{\scriptscriptstyle #1}}
\newcommand{\q}[1]{q^{\scriptscriptstyle #1}}
\renewcommand{\l}[1]{\lambda^{\scriptscriptstyle #1}}

\theoremstyle{plain}
\newtheorem{teo}{Teorem}
\newtheorem{theorem}[teo]{Theorem}
\newtheorem{proposition}[teo]{Proposition}

\theoremstyle{definition}
\newtheorem{definition}[teo]{Definition}

\theoremstyle{remark}
\newtheorem{remark}[teo]{Remark}

\begin{document}

\title[Spectral asymptotics of percolation Cayley graphs]{Spectral asymptotics of percolation Hamiltonians on amenable Cayley graphs}
\author[T.~Antunovi\'c]{Ton\'ci Antunovi\'c}

\author[I.~Veseli\'c]{Ivan Veseli\'c}
\address{Emmy-Noether-Programme of the Deutsche Forschungsgemeinschaft\vspace*{-0.3cm} }
\address{\& Fakult\"at f\"ur Mathematik,\, 09107\, TU\, Chemnitz, Germany   }

\urladdr{www.tu-chemnitz.de/mathematik/schroedinger/members.php}

\keywords{amenable groups, Cayley graphs, random graphs, percolation, random operators, spectral graph theory, phase transition}
\subjclass[2000]{05C25, 82B43, 05C80, 37A30, 35P15}

\begin{abstract}
In this paper we study spectral properties of adjacency and  Laplace operators 
on percolation subgraphs of Cayley graphs of amenable, finitely generated groups.
In particular we describe the asymptotic behaviour of the integrated density of states
(spectral distribution function) of these random Hamiltonians near the spectral minimum. 

The first part of the note discusses various aspects of the quantum percolation model, 
subsequently we formulate a series of new results, 
and finally we outline the strategy used to prove our main theorem.

\end{abstract}
\thanks{\jobname.tex; \today }
\thanks{\copyright 2007 by the authors. Faithful reproduction of this article,
        in its entirety, by any means is permitted for non-commercial purposes}
\maketitle

\section{Introduction}
The results on which we report rely on tools and ideas from various mathematical fields.
In the introduction we sketch the role they play in our study of percolation Hamiltonians.

\subsection{Integrated density of states} 	\label{ss-IDS}

For a large class of random operators which are ergodic with respect to a group of translations on the configuration space
an \textit{integrated density of states} (or \textit{spectral distribution function}) can be defined using 
an exhaustion of the whole space by subsets of finite volume. 
Since this fact relies on an ergodic theorem, it is not surprising that 
the underlying group needs to be amenable. Note that periodic operators are a special class of ergodic ones and
thus also posses a well defined integrated density of states (in the following abbreviated as IDS).
For ergodic random operators the spectrum is deterministic, i.e.~for any two realisations of the operator the 
spectrum (as a set) coincides almost surely. The same statement holds for the measure theoretic components of the spectrum.
However, there are other spectral quantities, like eigenvalues or eigenvectors, which are highly dependent on the randomness.
For random Schr\"odinger type operators on $L^2(\RR^d)$ and $\ell^2(\ZZ^d)$ all these fundamental results can be 
found e.g.~in the monographs \cite{CarmonaL-90,PasturF-92}. 
For the underlying original results of Pastur and Shubin for random, respectively almost-periodic operators see
e.g.~\cite{Pastur-71,Pastur-80} and \cite{Shubin-78,Shubin-79}.

Once the existence of the IDS is established, it is natural to ask whether, 
for specific  models, one can describe some of its characteristic features in more detail.
Among them the most prominent are: the continuity and discontinuity properties of the IDS, 
and the asymptotic behaviour near the spectral boundaries.
The interest in these questions has a motivation stemming from the quantum theory of solids.
A good unterstanding of the features of the IDS is an important step towards
the determination of the spectral types of the considered Hamilton operators. 
Of all the measure theoretic spectral types of random operators the pure point part
is so far understood best. 
In particular, for certain random Hamiltonians it is established
that there exists an energy interval with dense eigenvalues and no continuous spectral component, 
a phenomenon called \textit{localisation}. For detailed expositions of this subject see, for instance,
\cite{Stollmann-01,AizenmanSFH-01,GerminetK-03,AizenmanENSS-06}.
All the proofs of localisation known so far use as an essential ingredient 
estimates on the IDS or some closely related quantity.
Expository accounts of the IDS of various types of random Hamiltonians can be found among others in \cite{LeschkeMW-03}, 
\cite{KirschM-07} or \cite{Veselic-06b}.

\subsection{Lifshitz asymptotics}

The IDS of lower bounded periodic operators in Euclidean space exhibits typically a polynomial behaviour
at the minimum of the spectrum. This is for instance the case for 
the discrete Laplace operator on $\ell^2(\ZZ^d)$ and for Schr\"odinger operators on $L^2(\RR^d)$ with a periodic potential. 
Analogous results can be proven for uniformly elliptic divergence type operators with periodic coefficients 
using \cite{Kirsch-87b} and for quantum waveguides using e.g.~\cite{BirmanS-01} or \cite{KondejV-06b}.
In the mentioned cases the low energy behaviour 
of the IDS is characterised by the so called \emph{van Hove singularity}. This means that the IDS $N(\cdot)$ 
behaves for $0<E\ll 1$ asymptotically like $N(E)\sim E^{d/2}$, i.~e.~polynomially  with the \emph{van Hove exponent} equal to $d/2$ .

A random perturbation of the Hamilton operator changes drastically the low energy asymptotics. For many random models in 
$d$-dimensional Euclidean space it has been proven that the 
asymptotic behaviour of the IDS is exponential in the sense that $N(E) \sim \exp(-const (E-E_0)^{-d/2})$.
Here $E_0$ denotes the minimum of the spectrum of the random operator,
which is by ergodicity independent of the realisation almost surely.
The exponential behaviour  of the IDS  has been first deduced on physical grounds by {I.~M.~Lif\v sic} in \cite{Lifshitz-63,Lifshitz-64,Lifshitz-65} and is accordingly called 
\emph{Lifshitz asymptotics} or \emph{Lifshitz tail}.
The most precise bounds of this type have been obtained for random Schr\"odinger operators 
with a potential generated by impurities which are distributed randomly 
in space according to a Poisson process, see e.~g.~\cite{DonskerV-75c,Pastur-77,Nakao-77,Sznitman-98,KloppP-99}. 
Similar results hold for a discrete relative of this operator, namely the \emph{Anderson model} on $\ell^2(\ZZ^d)$, 
see e.g.~\cite{DonskerV-79,Antal-95,Klopp-00b,BiskupK-01a,Metzger-05,vanderHofstadKM-06}.
The reason why these models are amenable to a very precise analysis
is the applicability of Brownian motion, respectively random walk techniques and Feynman-Kac functionals.
For several other types of random operators the weaker form  
\be
    \label{e-LifTaildef}
    \lim_{E \searrow E_0} \frac{\log |\log (N(E)-N(E_0))|}{|\log (E-E_0)|} = \frac{d}{2}
\ee
of the exponential law has been established using rather simple estimates based on inequalities 
by Thirring \cite{Thirring-94} and Temple \cite{Temple-28}.
In particular, this method can be applied to show that the asymptotics \eqref{e-LifTaildef} 
holds for the Anderson model on $\ell^2(\ZZ^d)$ and its continuum counterpart, 
the \emph{alloy type model} on $L^2(\RR^d)$. 
This strategy of proof was pursued in the 1980s in several papers 
\cite{KirschM-83a,Mezincescu-85,Simon-85b,KirschS-86,Mezincescu-86,Mezincescu-87,Simon-87}
by Kirsch, Martinelli, Simon and Mezincescu.
These two models will serve as a point of reference in the present note, since 
they are closely related to percolation Hamiltonians. Moreover, the methods used to study the low-energy
asymptotics of percolation models are, at least in part, parallel to those of the last mentioned series of papers.
The main difference is that one needs to replace the Temple or Thirring bound by some inequality from spectral geometry.
We will call the right hand side of \eqref{e-LifTaildef} the \emph{Lifshitz exponent}.

Let us relate the asymptotic behaviour \eqref{e-LifTaildef} to two spectral features
mentioned already in \S \ref{ss-IDS}, namely spectral localisation and continuity of the IDS. 
If \eqref{e-LifTaildef} holds then the spectral edge $E_0$ is called a \emph{fluctuation boundary}.
This term stems from the fact, that if the system is restricted to a large, finite volume, 
eigenvalues close to $E_0$ correspond to very particular and rare realisations of the randomness. 
Thus the 
spectral edge on finite volumes is highly sensitive to fluctuations of the random configurations.
This feature is closely related to the phenomenon of localisation. Indeed,
for many models existence of pure point spectrum has been proven precisely 
in the energy regimes where the density of states is very sparse, 
see e.~g.~the characterisation given in \cite{GerminetK-04}.
The Lifshitz asymptotics implies that the IDS is extremely thin near the bottom of the spectrum, 
thus this energy region is a typical candidate for pure point spectrum.

One expects that the IDS is continuous at the fluctuation boundaries.
Indeed, the Lifshitz asymptotics implies continuity at the minimum $E_0$ of the spectrum
and moreover that the limes superior and inferior of difference quotients of any order 
vanish at $E_0$. 
This is in remarkable contrast to the dense set of discontinuities of the IDS 
exhibited by many percolation Hamiltonians, see e.g.~\cite{deGennesLM-59b,ChayesCFST-86,Veselic-05b}.

\subsection{Percolation Hamiltonians}

Hamiltonians on site percolation subgraphs of the lattice were introduced 
by P.-G.~de Gennes, P.~Lafore and J.~Millot in \cite{deGennesLM-59a,deGennesLM-59b} 
as quantum mechanical models for binary alloys.  
The resulting random operator bears a strong similarity to the Anderson model considered
by P.~W.~Anderson in \cite{Anderson-58}.  The difference is that in the site percolation
Hamiltonian of \cite{deGennesLM-59b} the random potential may assume only two values, 
namely zero and plus infinity.
The lattice sites where the potential equals infinity are deleted from the graph, and thus the 
Hamiltonian is restricted to the space of vertices where the potential vanishes. The resulting Hamiltonian 
may be understood as a single band approximation of an Anderson model whose potential is a family 
of Bernoulli distributed random variables.
More precisely, if one equips the Bernoulli-Anderson potential with a 
coupling constant and lets this one tend to infinity, the corresponding quadratic form converges
to the one of the quantum percolation model, cf.~Remarks \ref{rem: operator comparison} and \ref{rem: bands}. This limit is in the physical literature, 
cf.~e.g.~\cite{KirkpatrickE-72},
sometimes understood as the strong scattering limit of the tight binding model. It may be possible to estimate efficiently 
the convergence of the associated IDS using the techniques and estimates introduced in \cite{LenzMV}.

A series of papers in theoretical 
\cite{KirkpatrickE-72,ChayesCFST-86}
and computational physics, 
e.g~\cite{ShapirAH-82,KantelhardtB-97,KantelhardtB-98,KantelhardtB-98b,KantelhardtB-02}, 
analysed the spectral features of the quantum percolation model. A particular point of interest in the numerical
studies was to compare the spectral localisation of percolation Hamiltonians with the one of the Anderson model.

More recently there has been interest for such models in the mathematics community. The paper 	
\cite{BiskupK-01a}, which studies the intermittency behaviour of the parabolic Anderson model, establishes also
Lifshitz tails for the quantum site percolation model. Given the abovementioned relation between the 
Anderson and the quantum percolation Hamiltonian, it is not surprising that a detailed analysis of 
the former gives also results about the latter. Some related properties of random walks on percolation graphs 
have been analysed rigorously already in \cite{Antal-95}.

The existence of the IDS for a rather general class of site percolation Hamiltonians on graphs with a quasi transitive, 
free\footnote{In the references \cite{Veselic-05a,Veselic-05b} is was forgotten to spell out the assumption that the group acts freely on the graph.} 
amenable group action was established in \cite{Veselic-05a}.
This result relies on Lindenstrauss' pointwise ergodic theorem \cite{Lindenstrauss-01} for locally compact second countable amenable groups.
Likewise, the non-randomness of the spectrum and its components is valid also in this general setting, see \cite{Veselic-05b,LenzPV-02}.
These results hold in particular for the Anderson model and periodic Laplacians on such graphs.
In fact, these features of the percolation Hamiltonians and their proofs are quite analogous to those 
of the Anderson model on $\ell^2(\ZZ^d)$. However, there are some distinct features of the quantum percolation model
which distinguish it sharply from random Hamiltonians on the full lattice. The latter have a continuous IDS 
\cite{DelyonS-84,CraigS-83a}, while percolation Hamiltonians have a dense set of discontinuities, which 
survives even if one restricts the operator to infinite percolation clusters \cite{KirkpatrickE-72,ChayesCFST-86,Veselic-05b}.
In \cite{KirschM-06} Kirsch \& M\"uller analysed basic spectral features of bond percolation Laplacians on the lattice
and moreover carried out a thorough study of the low-lying spectrum of these operators in the non-percolating regime.
These results were complemented by M\"uller \& Stollmann in a paper \cite{MuellerS} where the percolating regime is studied.
The present note continues the analysis \cite{Veselic-05b} of site percolation Hamiltonians on general graphs and at the same time extends
the results of \cite{KirschM-06} to bond percolation models on amenable Cayley graphs. In particular, 
part of our proofs relies on ideas introduced in \cite{KirschM-06} and extends them to a more general geometric setting.

It is not surprising that the analysis of (certain) spectral properties of percolation Hamiltonians
relies on the proper understanding of the underlying 'classical' percolation problem.
An exposition of this theory applied to independent bond percolation on the 
lattice can be found in Grimmett's book \cite{Grimmett-99}. Some other aspects are covered in 
the monograph \cite{Kesten-82} by Kesten.
In the more recent literature quite a body of work is devoted to percolation processes on more general graphs than $\ZZ^d$,
see for instance \cite{BenjaminiS-96a,BenjaminiLPS-99b,LyonsPeres,ProcacciB-04}. Still, much more is known for percolation 
on lattices than on general graphs.
For the purposes of the results presented in this paper, 
we extended the theorems in \cite{AizenmanN-84,AizenmanB-87,Menshikov-86} 
on the sharpness of the phase transition and the exponential decay of the cluster size in the subcritical regime to 
quasi transitive graphs, see \cite{AntunovicV-a}.

Since in the subcritical regime no infinite percolation cluster exists almost surely,
the entire spectrum of the corresponding Hamiltonian consists of eigenvalues.
In the supercritical regime an infinite cluster exists and thus 
the Hamiltonian may have continuous spectrum. It is however conjectured \cite{ShapirAH-82} that 
for values of the percolation parameter just above the critical point the 
spectrum will still have no continuous component, although an infinite cluster exists.
Thus there may be a second `quantum' critical value of the percolation parameter,
strictly larger than the classical critical value.
Let us note that there is no rigorous proof of Anderson localisation for percolation Hamiltonians (in the percolating regime).
The reason is that this model has Bernoulli distributed randomness, as it is the case for the Bernoulli-Anderson Hamiltonian.
Due to the singular nature of the distribution of the randomness the known proofs of localisation do not apply.

\subsection{Spectral graph theory and geometric $L^2$-invariants}
\label{ss-sgt}
To make sense of the term `low energy asymptotics' one has to know 
where the minimum of the spectrum lies. In the case of Cayley graphs of amenable groups
it is known from Kesten's theorem \cite{Kesten-59a} that
the  bottom of the spectrum of the Laplacian is equal to zero.
If the graph under consideration happens to be bipartite the spectrum of the adjacency operator 
is symmetric with respect to the origin, see e.g.~\cite{CoulsonR-40,Simon-85b}. 
This allows one to translate 
results about the lower spectral boundary of adjacency and Laplace operators 
to statements on the upper boundary.

The low energy behaviour of the IDS is determined by the first non-zero eigenvalue 
of the Hamiltonian on percolation clusters with `optimal shape'.
The minimising configurations are determined by inequalities from spectral graph theory 
like the isoperimetric, the Cheeger and the Faber \& Krahn inequality, 
cf.~\cite{CoulhonSC-93,Chung-97,ColindeVerdiere-98,ChungGY-00}.
They relate the lowest eigenvalue on a finite subgraph to its volume and the 
volume growth behaviour of the Cayley graph. The growth of balls on Cayley graphs, in turn, 
can be classified using a result of Bass \cite{Bass-72} and Gromov's theorem \cite{Gromov-81}. We use
an improved version of the latter due to van den Dries \& Wilkie \cite{vandenDriesW-84}.
Let us mention that even if one is interested in purely geometric 
properties of the automorphism group of a graph, like amenability or unimodularity,
percolation may be a tool of choice, see for instance \cite{BenjaminiLPS-99b}.
There several  equivalences are established, each connecting a geometric property of a graph
with a probabilistic property of an associated percolation process.

For several types of Hamiltonians on the lattice the Lifshitz exponent of the 
random operator equals the exponent of the van Hove singularity 
of its periodic counterpart. It turns out that the same holds for 
certain percolation models on Cayley graphs. To formulate this more properly one needs to characterise 
the high energy behaviour of the adjacency operator (resp.~low energy behaviour of the Laplacian) on the full Cayley graph. 
This can be done by relating it to known results on random walks on groups, cf.~\cite{Woess-00},
or on geometric $L^2$-invariants, cf.~e.g.~\cite{Lueck-02}. In fact, it turns out that the van Hove exponent of a
Cayley graph equals the first \emph{Novikov-Shubin invariant} of the Laplacian on the graph.
These invariants have been introduced in  \cite{NovikovSh-86,NovikovSh-86a} and studied e.g.~in \cite{GromovS-91}.
For various analogies between geometric $L^2$-invariants and properties of the IDS see
the workshop report \cite{DodziukLPSV-06}. The abovementioned equality of the Lifshitz and 
the van Hove exponent is encountered also in other contexts, see for instance Klopp's analysis
\cite{Klopp-99} of Lifshitz tails at internal spectral gap edges. 

\subsection{The main result}
Let us loosely state the main result of this note. 
Consider a Cayley graph of a finitely generated amenable group.
Assume that the volume of balls of radius $n$ in the graph 
behaves like $n^d$ with the convention that $d=\infty$ corresponds to superpolynomial growth.
Each deleted site (respectively bond) induces a new boundary in the graph, at which we may impose 
a certain type of boundary condition giving rise to different Laplace operators. More precisely,
since we are dealing with bounded operators, the boundary condition is not a restriction 
on the domain of the operator, but rather an additional boundary term.
As in \cite{KirschM-06} we consider Dirichlet, Neumann, and adjacency (or pseudo-Dirichlet) 
percolation Laplacians. 

For the adjacency and Dirichlet Laplacian the low energy asymptotics of the IDS is given by 
\be
    \lim_{E \searrow 0} \frac{\log |\log N(E)|}{|\log E|} = \frac{d}{2}
\ee
Thus we have a Lifshitz type behaviour and zero is a fluctuation boundary of the spectrum.
The Lifshitz exponent coincides with the van Hove exponent of the underlying full Cayley graph.
For the Neumann Laplacian the low energy asymptotics of the IDS is given by 
\be
    \lim_{E \searrow 0} \frac{\log |\log (N(E)-N(0))|}{|\log E|} = \frac{1}{2}
\ee
where $N(0)$ is a non-zero value corresponding to the number of open clusters per vertex in the percolation graph.

\section{Model and results \label{ss-results}}

We present several results which we have obtained for the IDS of Hamiltonians on full Cayley graphs
and on graphs diluted by a percolation process. For the analysis of the latter operators it was necessary to 
establish certain properties of the `classical percolation model' on  Cayley graphs, which are detailed below.

\subsection{The Laplace Hamiltonian on a full Cayley graph}

We consider Cayley graphs of finitely generated, amenable groups and the corresponding Laplace operator.
Its IDS exhibits a van Hove asymptotics at the lower spectral edge whose exponent is determined by the 
volume growth behaviour of the group.

To formulate this more precisely we explain the geometric setting in detail.
Let $\Gamma$ be a finitely generated group. Each choice of a finite, symmetric 
set of generators $S$ not containing the unit element $\iota$ of $\Gamma$ gives rise to Cayley graph.
Note that this graph is regular (i.e.~all vertices have the same degree) 
with degree equal to the number of elements in $S$ and that the vertex set of $G$ can be identified with $\Gamma$.
Using the distance function on the Cayley graph we define the ball $B(n)$ 
of radius $n$ around the unit element $\iota$ in $\Gamma$ and set $V(n):=|B(n)|$. 
For a positive integer $d$ we use the notation $V(n) \sim n^{d}$ 
to signify that there exist constants $0 <a,b < \infty$ such that $a\, n^{d} \leq V(n) \leq b\, n^{d}$.

From now on we tacitly assume that all considered groups are finitely generated.
The growth of all such groups can be classified using deep
results of Bass \cite{Bass-72}, Gromov \cite{Gromov-81} and van den
Dries \& Wilkie \cite{vandenDriesW-84}. 

\begin{theorem}\label{thm: growth of groups} 
Let $G$ be a Cayley graph of a finitely generated group. Then exactly one of the following is true:
\begin{enumerate}[(a)]
\item $G$ has polynomial growth, i.e.~$V(n) \sim n^{d}$ holds for some $d \in \NN$,
\item $G$ has superpolynomial
growth, i.e.~for every $d\in\NN$ and every $b\in\RR$ there exists only finitely many integers $n$
such that $V(n) \leq  b\, n^{d}$.
\end{enumerate}
The growth behaviour (in particular the exponent $d$) is independent of the chosen set of generators $S$.
\end{theorem}

More precisely, the theorem follows from the following results.
Bass has shown that every nilpotent group satisfies $V(n) \sim n^{d}$ for some $d \in \NN$,
sharpening earlier upper and lower bounds of Wolf \cite{Wolf-68}.
Gromov's result in \cite{Gromov-81} is that every group of polynomially bounded 
volume growth is virtually nilpotent. Finally, van den Dries \& Wilkie have shown 
that the conclusion of Gromov's theorem still holds, if one requires 
the polynomial bound $V(n_k) \le b \, n_k^d$ merely along
a strictly increasing sequence $(n_k)_{k\in\NN}$ of integers $n_k \in \NN$.
We note that Pansu \cite{Pansu-83} has shown that $c:=\lim_{n\to\infty} V(n)/n^d$ exists.

Let us now define the Laplacian on $G$. For future reference we consider a more general
situation than needed at this stage. Let $G=(V,E)$ be a connected regular graph of degree $k$
with vertex set $V$ and edge set $E$, and $G'=(V',E')$ an arbitrary subgraph of $G$. Note that $G'$ 
is in general not regular. We denote the degree of the vertex $x\in V'$ in $G'$ by $d_{G'}(x)$.
If two vertices $x,y\in V'$ are adjacent in the subgraph $G'$ we write  $y \sim_{G'} x$.

\begin{definition}\label{def: laplace operator}
For $G$ and $G'$ as above we define the following operators on $\ell^{2}(V')$.
\begin{enumerate}[(a)]
\item The identity operator on $G'$ is denoted by $\Id_{V'}$. 
If there is no danger of confusion we drop the subscript $V'$.
\item The \emph{degree operator} of $G'$ acts on $\varphi \in\ell^{2}(V')$ according to
 \[
 [D(G')\varphi](x) := d_{G'}(x)\varphi(x).
 \]
\item The \emph{adjacency operator} on $G'$ is defined as 
 \[
 \displaystyle [A(G')\varphi](x) := \sum_{\genfrac{}{}{0pt}{2}{y \in V'}{y \sim_{G'} x}}\varphi(y).
 \] 
\item The \emph{adjacency Laplacian}  on $G'$ is defined as
 \[
 \D{A}(G') := k \Id_{V'} - A(G').
 \]
\item In the special case $G'=G$, the \emph{Laplacian} or \emph{free Hamiltonian} on $G$ is defined as
 \[
 \Delta(G) := k \Id_{V} - A(G).
 \] 
(Thus it coincides with the adjacency Laplacian on $G$.)
\end{enumerate}
\end{definition}

There are several different names used for the adjacency Laplacian 
$\D{A}(G') $ in the literature.  Our terminology is motivated by the fact that
up to a multiplicative and an additive constant it is equal to the adjacency operator 
on the subgraph $G'$. This will be different for 
operators with an additional Dirichlet or Neumann boundary term introduced in Subsection \ref{ss-percolation Hamiltonians}.
For induced subgraphs $G'$ of $G$ one can consider the 
\emph{restriction of the Laplacian} which is defined as 
$  \D{P}(G') := P_{V'} \Delta(G) P_{V'}^*$. 
Here $P_{V'}\colon \ell^2(V)\to \ell^2(V'), P_{V'}\varphi(x):=\chi_{V'}(x)\varphi(x)$ 
denotes the orthogonal projection on $V'$.
It turns out that $\D{P}(G') =k  \Id_{V'} - A(G')= \D{A}(G')$.

The IDS of the {free Hamiltonian} on a Cayley graph $G$ of an amenable group $\Gamma$ 
can be defined as $N_0(E):= \la\chi_{]-\infty,E]}(\Delta(G))\delta_{\iota},\delta_{\iota}\ra$.
Here the function $\delta_{\iota}$ has value $1$ at $\iota$ and $0$ everywhere else. 
It is possible to construct $N_0$ via an exhaustion procedure, see e.g.~\cite{Lueck-94c,DodziukM-98,DodziukLMSY-03}.

The next Theorem characterises the asymptotic behaviour of the IDS at the spectral bottom.
For groups of polynomial growth it exhibits a van Hove singularity. 
For groups of superpolynomial growth one encounters a type of generalised 
van Hove asymptotics with the van Hove exponent equal to infinity.

\begin{theorem}\label{thm: full graph}
Let $\Gamma$ be a finitely generated,  amenable group, $\Delta(G)$ the Laplace operator on a Cayley graph of $\Gamma$
and $N_0$ the associated IDS. If $\Gamma$ has polynomial growth of order $d$ then
\begin{align}
\label{e-vanHove}
    \lim_{E \searrow 0} \frac{\log N_0(E)}{\log  E} &= \frac{d}{2}.
\intertext{and if $\Gamma$ has superpolynomial growth then} 
\label{e-vanHove infinite}
    \lim_{E \searrow 0} \frac{\log N_0(E)}{\log  E} &= \infty.
\end{align}
\end{theorem}

\subsection{Classical percolation\label{ss-classical-percolation}}

Next we briefly introduce percolation on graphs. More precisely, we will consider 
independent site percolation, as well as independent bond percolation on \emph{quasi transitive} graphs.
These are graphs whose vertex set decomposes into finitely many equivalence classes under 
the action of the automorphism group. Such graphs have uniformly bounded vertex degree.
If there is only one orbit, the graph is called \emph{transitive}. 
Cayley graphs are particular examples of transitive graphs.

Let $G=(V,E)$ be an infinite, connected, quasi transitive graph and $p$ some fixed real number between $0$ and $1$. 
For every site $x \in V$ we say it is \emph{open} with probability $p$ and \emph{closed} with probability $1-p$ 
independently of all other sites. This is the most simple site percolation model.
More formally, we consider for every vertex $x$ the probability space $\Omega_{x}:=\left\{0,1\right\}$, 
with the sigma algebra consisting of the power set $\cP(\Omega_{x})$ of $\Omega_{x}$ and the probability measure 
$\mathbb{P}_{x}$ on $(\Omega_{x}, \cP(\Omega_{x}))$ defined by $\mathbb{P}_{x}(0)=1-p$ and
$\mathbb{P}_{x}(1)=p$. The probability space 
$(\Omega_V,\mathcal{F}_V,\PP_V)$ associated to site percolation is defined as the product 
$\prod_{x \in V}(\Omega_{x},\cP(\Omega_{x}),\mathbb{P}_{x})$.
The coordinate map $\Omega_V\ni\omega \mapsto \omega_x$ defines the \emph{site percolation process} $\Omega_V\times V\to \{0,1\}$.
We call $V(\omega):=\{x \in V\mid \omega_x =1\}$ the set of open or active sites. The induced subgraph of $G$
with vertex set $V(\omega)$ is denoted by $G(\omega)$ and called the \emph{percolation subgraph} 
in the configuration $\omega$. The connected components of $G(\omega)$ are called \emph{clusters}.

The bond percolation process is defined completely analogously.
The bond percolation probability space is 
$(\Omega_E,\mathcal{F}_E,\mathbb{P}_E) := \prod_{e \in E}(\Omega_{e},\cP(\Omega_{e}),\mathbb{P}_{e})$,
where for every edge $e \in E$ we have $\Omega_{e}:=\left\{0,1\right\}$, with power set $\cP(\Omega_{e})$ 
and probability measure $\mathbb{P}_{e}$ defined by $\mathbb{P}_{e}(0)=1-p$ and $\mathbb{P}_{e}(1)=p$. 
For a given configuration $\omega \in\Omega_E$ we define the \emph{percolation subgraph}  
$G(\omega)$ as the graph whose edge set $E(\omega)$ is the set of all $e\in E$ with $\omega_e =1$ and whose
vertex set $V(\omega)$ consist of all vertices in $V$ which are incident to an element of $E(\omega)$.
Connected components of $G(\omega)$ are called clusters. We will denote the expectation with respect to 
either $\PP_V$ or $\PP_E$ by $\EE \{  \dots\}$.

The most basic result in percolation theory is that for both the site and the bond model
there exists a critical parameter $0<p_{c}\leq 1$ such that the following statement holds:
\begin{itemize}
\item if $p<p_{c}$ there is no infinite cluster almost surely
(\emph{subcritical phase}),
\item if $p>p_{c}$ there is an infinite cluster almost surely (\emph{supercritical phase}).
\end{itemize}
Of course, the value of $p_c$ depends on the graph and on the type of percolation process considered.

It turns out that more can be said about the size of clusters in the subcritical phase.
Let $o\in V$ be an arbitrary, but fixed vertex in a quasi transitive graph $G$.
Consider site or bond percolation on $G$ and denote by $C_o(\omega)$ the cluster containing $o$
and by $|C_o(\omega)|$ the number of vertices in $C_o(\omega)$.
\begin{theorem}\label{thm: exponential decay}

For any $p<p_{c}$ there exist a constant $\tau_{p}>0$ 
such that $\PP( |C_{o}(\omega)| \geq n) \leq e^{-\tau_{p}n}$ for all $n\in \NN$. In particular, 
the expected number of vertices in $C_{o}(\omega)$ is finite for all $ p< p_c$.
\end{theorem}

This extends earlier results of Aizenman \& Newman \cite{AizenmanN-84}, Menshikov \cite{Menshikov-86} 
(see also \cite{MenshikovMS-86}), and Aizenman \& Barsky \cite{AizenmanB-87}.
The generalisation to bond and site percolation on quasi transitive graphs is given in \cite{AntunovicV-a} 
and relies on the differential inequalities for order parameters established in \cite{AizenmanB-87}.
In our later applications, we will need this result only for Cayley graphs of 
amenable finitely generated groups.

\subsection{Percolation Hamiltonians on Cayley graphs}
\label{ss-percolation Hamiltonians}

We introduce percolation Hamiltonians with and without boundary terms and 
establish the existence of a self-averaging IDS. This enables us to state our main result about 
the Lifshitz asymptotics and to compare it with the van Hove singularities of the 
Hamiltonian on the full Cayley graph.
\smallskip

Consider an arbitrary subgraph $G'=(V',E')$ of an infinite, $k$-regular Cayley graph $G=(V,E)$.
The multiplication operator 
\[
W^{b.c.}(G')= k \Id_{V'} -D(G')
\]
is non negative and has support in the interior vertex boundary of $G'$. 
Thus $\pm W^{b.c.}(G')$  can be understood as a potential due to the repulsion/attraction of the boundary.
When added to the adjacency Laplacian $\Delta^A$ it gives rise to Dirichlet and Neumann boundary condition.
Here we follow the nomenclature of \cite{Simon-85b,Mezincescu-86,KirschM-06}. 

\begin{definition}\label{def: boundary conditions}
We define the following operators on $\ell^{2}(V')$
\begin{enumerate}[(a)]
\item The \emph{Dirichlet Laplacian} is defined as 
 \[
 \D{D}(G') := \Delta^A +W^{b.c.}(G')= 2k \Id_{V'} - D(G') - A(G').
 \]
\item The \emph{Neumann Laplacian} is defined as 
 \[
 \D{N}(G') := \Delta^A -W^{b.c.}(G')= D(G') - A(G').
 \]
\end{enumerate}
\end{definition}
The Laplacian $\Delta(G)$ on the full graph will be abbreviated by $\Delta$. Note that in this case 
all three versions of the Laplacian coincide since there is no boundary.
Let us collect certain basic properties of the adjacency, Dirichlet and Neumann Laplacian.

\begin{remark}\label{rem: properties}
\newcounter{Lcount}
\begin{list}{(\alph{Lcount})}{\itemindent -2,75em} \usecounter{Lcount}
\item 
All operators introduced in Definitions \ref{def: laplace operator} and \ref{def: boundary conditions} 
are bounded and self-adjoint. 
\item  Since $W^{b.c.}(G') \ge 0$ we have the following inequalities in the sense of quadratic forms:
\begin{equation}\label{eq:relations}
\D{N}(G') \leq \D{A}(G') \leq \D{D}(G') 
\end{equation}
This is complemented by $0 \leq \D{N}(G') $ and $\D{D}(G') \leq 2k \Id$, 
which can be show by direct calculation in the same manner as in \cite{KirschM-06}.
Later we will see, that there is a more complete chain of inequalities between five operators.
\item The operators $\D{A}(G')$ and $\D{D}(G')$ are injective. 
The Neumann Laplacian $\D{N}(G')$ is injective if and only if $G'$ has no finite components. 
In particular, the dimension of the kernel of $\D{N}(G')$ is the number of the finite components of $G'$.
Let us note that this relation between the number of zero Neumann eigenvalues 
and the number of clusters makes it possible to refine the analysis of the large time behaviour 
of random walks on finite percolation clusters, see e.g.~\cite{Sobieczky}.
\item As we said in \S \ref{ss-sgt}, it is a classical result \cite{Kesten-59a} that zero is the lower spectral edge 
of $\Delta$ if and only if $G$ is an amenable graph. Finding lower bounds in the nonamenable case 
is a more difficult task, see for instance \cite{Zuk-97,BartholdiCCSdlH-97,Nagnibeda-97}. 
\item A graph $G$ is  called \emph{bipartite} if there is a partition $V_1, V_2$ of the vertex set 
such that there is no edge in $G$ which joins two elements of $V_i$ for $i=1,2$. For bipartite graphs 
it is useful to consider conjugation with the operator $\cU$ on $\ell^{2}(V)$ which is given by the
multiplication with the function $\chi_{V_1}-\chi_{V_2}$.
This operator is unitary, selfadjoint and an involution. 
We will denote the restriction of $\cU$ to some subset $\ell^2(V')\subset \ell^2(V)$ again by the
same symbol. The conjugation with $\cU$ 
relates the upper spectral edge to the lower one, see e.g.~\cite{Simon-85b,Mezincescu-86,KirschM-06}.
For any subgraph $G'$ we have $D(G')=\cU D(G')\,\cU$ and $A(G')=-\cU A(G')\,\cU$. This implies the relations
\begin{equation}\label{eq: bipartite relations}
\D{A}(G') = 2k\Id - \cU\D{A}(G')\,\cU\quad \text{ and } \quad \D{N(D)}(G') = 2k\Id - \cU\D{D(N)}(G')\,\cU.
\end{equation}
Consequently the spectrum of $\D{A}(G')$ is symmetric with respect to $k$, while the spectrum of
$\D{N}(G')$ is a set obtained from the spectrum of $\D{D}(G')$ by the symmetry with respect to $k$. 
It follows that in the special case of amenable bipartite graphs the upper spectral edge of $\Delta$ is equal to $2k$.
Bipartiteness is not only sufficient, but also necessary for the mentioned symmetry. 
Claim 4.5 from \cite{Oguni-07} shows that on amenable
non-bipartite graphs $2k$ is not in the spectrum of $\Delta$, thus the symmetry fails.
\end{list}\end{remark}
\smallskip

Now we introduce percolation Laplacians associated to the configuration $\omega$
as bounded, selfadjoint operators  on $\ell^2(V(\omega))$.

\begin{definition}\label{def: percolation laplacians}
Let $(\Omega, \mathcal{F}, \mathbb{P})$ be a probability space corresponding 
either to site or bond percolation on $G$. For $\omega \in \Omega_V$ (resp.~$\omega \in \Omega_E$)
the operator 
\[
\D{\#}_\omega:=\D{\#}(G(\omega))\colon \ell^{2}(V(\omega))\to \ell^{2}(V(\omega)), 
\quad  {\scriptstyle\#}\in \left\{A,D,N\right\}
\]
is called \emph{adjacency, resp.~Dirichlet, resp.~Neumann percolation Laplacian}.
\end{definition}
For brevity we introduce the notation $A_\omega := A(G(\omega))$ and $W_\omega^{b.c.}:=W^{b.c.}(G(\omega))$.
Note that in the site percolation model we have 
$\displaystyle d_{G(\omega)}(x):=\omega_x \sum_{y \sim x} \omega_y$,
while for bond percolation 
 \[ 
 d_{G(\omega)}(x):=\chi_{V(\omega)}(x) \sum_{e \sim x} \omega_e
 \]
where $e \sim x$ signifies that in the graph $G$ the edge $e$ is incident to $x$.
Note that if the distance of $x,y \in V$ is greater than one, the random variables 
$d_{G(\omega)}(x)$ and $d_{G(\omega)}(y)$ are independent.
These facts imply that in both models the stochastic field $(\omega,x)\mapsto W_\omega^{b.c.}(x) $
gives rise to a random potential which is stationary and ergodic with respect to the action of the group $\Gamma$.
\smallskip

Due to the structure of the underlying percolation process the considered random operators 
satisfy an equivariance relation which we explain next. 
For a group $\Gamma$ and an associated Cayley graph $G$ there is a natural action of $\Gamma$ 
on $G$ by multiplication from the left. 
This gives rise to a $\Gamma$-action both on the corresponding percolation probability space 
and on the bounded operators on the $\ell^2$ space over the graph.
The operation on the site percolation probability space $(\Omega_V,\PP_V)$
by measure preserving transformations is given by $(\tau_{\gamma}(\omega))_{x} := \omega_{\gamma^{-1}x}$, 
for $\gamma\in \Gamma$ and $x \in V$. Note that $\gamma$ and $x$ are from the same set, but the
first one is considered as a group element, while the second as a vertex. The family of
transformations $(\tau_{\gamma})_{\gamma \in \Gamma}$ acts ergodically on $(\Omega_V,\PP_V)$.
In the same way there is an ergodic action of measure preserving transformations indexed 
by $\gamma \in \Gamma$, which we again denote by $\tau_\gamma$, on the 
bond percolation probability space $(\Omega_E,\PP_E)$. Here the transformations are given by 
$(\tau_{\gamma}(\omega))_{e} := \omega_{\gamma^{-1}e}$, where $\gamma^{-1}[y,z] = [\gamma^{-1}y,\gamma^{-1}z]$
and $[y,z]$ denotes the edge joining the vertices $y$ and $z$.

The group $\Gamma$ acts by unitary translation operators $(U_{\gamma}\varphi)(y):=\varphi(\gamma^{-1}y)$ on $\ell^{2}(V)$. 
If we restrict these operators to the $\ell^2$ spaces over the active sites they act consistently 
with the shift on the probability space, more precisely 
$U\colon \ell^2(V(\omega)) \to \ell^2(\gamma V(\omega))=\ell^2(V(\tau_\gamma\omega))$. 
Due to the transformation behaviour of the adjacency operator $A_\omega$
and the boundary term  potential $W_\omega^{b.c.}$ one has the equivariance relation 
$\D{\#}(\tau_{\gamma}\omega) = U_{\gamma}\D{\#}(\omega)U_{\gamma}^{-1}$ for ${\scriptstyle\#}\in\{A,D,N\}$. 
Since $(\tau_{\gamma})_{\gamma \in \Gamma}$ acts ergodically on $\Omega$, 
$(\D{\#}_\omega)_\omega$ falls into the class of random operators studied in \cite{LenzPV-02}.

This yields the non-randomness of the spectrum and spectral components. 
In the following let us denote by $\sigma$ the spectrum, 
and by $\sigma_{disc}$, $\sigma_{ess}$, $\sigma_{pp}$, $\sigma_{sc}$, $\sigma_{c}$  and $\sigma_{ac}$, 
the discrete, the essential, the pure point, the singular continuous, the continuous, 
and the absolutely continuous component of the spectrum respectively.
Let us recall that the pure point spectrum of an operator is the \emph{closure} of the set of its eigenvalues.
We use the notation $\sigma_{fin}(H)$ for the set of eigenvalues of $H$ which posses an eigenfunction with compact, 
i.e.~finite, support.

The following theorem holds for site and bond percolation Hamiltonians and for all values of $p \in[0,1]$.

\begin{theorem}\label{thm: nonrandom spectrum}
Let $\D{\#}_\omega$  be one of the percolation Laplacians defined above. 
Then there exists for each $\bullet \in \{disc,ess,pp,sc,c,ac,fin\}$ a subset 
$\Si{\#}_{\bullet}\subset\mathbb{R}$ and a subset $\Omega' \subset \Omega_V$ (resp.~$\Omega' \subset \Omega_E$) 
of full measure, such that $\sigma_{\bullet}(\D{\#}_\omega)=\Si{\#}_{\bullet}$ 
holds for every $\omega \in \Omega'$ and for all $\bullet \in \{disc,ess,pp,sc,c,ac,fin\}$. 
In particular $\sigma(\D{\#}_\omega)=\Si{\#}$ for all $\omega \in \Omega'$.
\end{theorem}

This has been proven for site percolation models in \cite{Veselic-05a,Veselic-05b,LenzPV-02}, 
and holds with the same proofs for bond percolation. 
For similar results in the case of bond percolation on the lattice see also \cite{KirschM-06}.

\begin{remark}\label{rem: spectrum}
(a) All the results stated so far hold for arbitrary quasi transitive graphs. In
particular it is not necessary to assume that the graph is amenable or infinite. 
For infinite graphs the discrete spectrum $\Sigma_{disc}$ is empty.

(b) The percolation Laplacians $\D{\#}_\omega$ are defined on $\ell^2(V(\omega))$. Thus, technically 
$(\D{\#}_\omega)_{\omega\in \Omega}$ is associated to a direct integral operator with non-constant fibre.
One can modify the operator in such a way that the fibres become constant and that the 
spectrum changes in a controlled way.
First note that for a fixed configuration $\omega$ and a fixed cluster $C(\omega)$
any percolation Hamiltonian leaves $\ell^2(C(\omega))$ invariant. Consequently, the full operator
decomposes according to the clusters:
\be \label{e:decomposition}
\D{\#}_\omega= \bigoplus_{C(\omega) \text{ cluster of } G(\omega)} \Delta^{\scriptstyle\#}(C(\omega)))
\ee
For any constant $K$ we may add to \eqref{e:decomposition} the operator 
$\bigoplus_{x \in V \setminus V(\omega)} K \Id_x = K \Id_{V \setminus V(\omega)}$. 
Denote the resulting sum by $\widetilde{\D{\#}_\omega}$. It acts on $\ell^2(V)$ 
and leaves the subspaces  $\ell^2(V(\omega))$  and $\ell^2(V\setminus V(\omega))$ invariant. In particular,
$\widetilde{\D{\#}_\omega}$ can be written in the same way as \eqref{e:decomposition} 
with the only difference that the direct sum extends over all $\tilde C(\omega)$ which are either a cluster
in $G(\omega)$ or a vertex in $V\setminus V(\omega)$. If the configuration $\omega$ is such that 
$V(\omega)=V$ then $\widetilde{\D{\#}_\omega}$ coincides with $\D{\#}_\omega$. 
Note however that these configurations form a set of measure zero, both in the site and the bond percolation model.
In all other cases we have 
\begin{align*}
\sigma(\widetilde{\D{\#}_\omega})&=\sigma(\D{\#}_\omega) \cup \{K\},
\\
\sigma_{pp}(\widetilde{\D{\#}_\omega})&=\sigma_{pp}(\D{\#}_\omega) \cup \{K\} \quad \text{ and }
\\
\sigma_{c}(\widetilde{\D{\#}_\omega})&=\sigma_{c}(\D{\#}_\omega).
\end{align*}

(c) In the subcritical regime there are no infinite clusters. The decomposition \eqref{e:decomposition} implies that
the operator is a direct sum of finite dimensional operators, thus the spectrum consists entirely of eigenvalues in 
$\sigma_{fin}(\D{\#}_\omega)$. In particular, $\Sigma_{sc} = \Sigma_{ac} =\emptyset$ and  $\Sigma_{pp} = \Sigma_{fin}$.

(d) We are not able to calculate the deterministic spectrum $\Si{\#}$ but only to give partial information about it.
Since we can find arbitrarily large clusters in $G(\omega)$ almost surely, strong resolvent convergence gives that 
$\Si{\#} \supset \sigma(\Delta)$, see \cite{KirschM-06} for details. 
Together with $\sigma(\D{\#}_\omega)\subset [0,2k]$ for all $\omega$ 
this implies that in the case of amenable Cayley graphs zero is the lower spectral edge of the operators 
$\D{\#}_\omega$ almost surely  for any ${\scriptstyle\#}\in \left\{A,D,N\right\}$.
Similar inclusions for the deterministic spectrum of Anderson models have been given in \cite{KunzS-80}.
If the Cayley graph is bipartite and the spectrum of the free Laplacian $\Delta$ has no gaps we have
$\Si{\#}= [0,2k]$. If $\sigma(\Delta)$ has several components, it might be that 
$\Si{\#}$ contains values in a spectral gap of $\Delta$. 
This is related to the phenomenon called \emph{spectral pollution} 
encountered when using strong convergence to
approximate spectral values of the limiting operator,
see for instance \cite{DaviesP-04,LevitinSh-04} and the references therein.
\end{remark}

Our next aim is to introduce the IDS and state its main properties. 
An abstract IDS may be defined without an exhaustion procedure by setting
\[
\N{\#}(E) = \EE\{\la\chi_{\left]-\infty,E\right]}(\D{\#}_\omega)\delta_{\iota},\delta_{\iota}\ra\}
\] 
see \cite{LenzPV-02}. However, the physically interesting situation is when 
the IDS is \emph{selfaveraging}, i.e.~when spatial averages of the normalised eigenvalue counting functions
 converge to the expectation with respect to the randomness. 
To be able to show that the IDS has indeed this property we have to require that the group $\Gamma$ is amenable.
Each such group has a sequence of finite, non-empty  subsets $(I_{j})_{j}$ which satisfies
for every finite set $F \subset \Gamma$ the property 
 \[
 \lim_{j \to \infty} \frac{|I_{j}\,\triangle \,F \cdot I_{j}|}{|I_{j}|} = 0.
 \] 
Such an $(I_{j})_{j}$ is called a \emph{F\o lner sequence}.  
Here $\triangle$ denotes the symmetric difference of two sets.
Various properties of such sequences are discussed in \cite{Adachi-93},
and their role in the construction of the IDS of random operators in \cite{PeyerimhoffV-02} 
and \cite[\S 2.3]{Veselic-04a}. 
Each F\o lner sequence has a tempered subsequence $(I_{j_{k}})_{k}$, meaning that there exists a constant $C$ 
such that $|I_{j_{k}}I_{j_{k-1}}^{-1}| \leq C |I_{j_{k}}|$.
We may consider $I_{j_k}$ as a subset of the vertices in the Cayley graph. 
Denote by $P_k$ the projection on the vertex set $I_{j_k}\cap V(\omega)$ and by 
$\D{\#,k}_\omega$ the restricted operator $P_k \D{\#}_\omega P_k^*$
for any ${\scriptstyle\#}\in \left\{A,D,N\right\}$. Thus we obtain a selfadjoint operator on a 
finite dimensional space which has a finite number of real eigenvalues. 
Hence the trace of the spectral projection $\Tr \big[\chi_{]-\infty,E]} (\D{\#,k}_\omega)\big]$ is finite.
We define the eigenvalue counting distribution function as 
 \[
 \N{\#,k}_\omega(E):=\frac{1}{|I_{k}|} \,\Tr \big[\chi_{]-\infty,E]} (\D{\#,k}_\omega)\big].
 \]
The next theorem holds for site and bond percolation Hamiltonians and for all values of $p \in[0,1]$.

\begin{theorem}\label{thm: IDS}
Let $G$ be a Cayley graph of an amenable group $\Gamma$ and $(I_{k})_{k}$ a tempered F\o lner sequence. 
There exists a set $\Omega' \subset \Omega_V$ (resp.~$\Omega_E$) 
of full measure such that for every $\omega \in\Omega'$ and every $E\in\RR$
which is a continuity point of $\N{\#}(\cdot)$  we have
\[
\lim_{k\to\infty} \N{\#,k}_\omega(E)= \N{\#}(E).
\] 
The support of the measure associated to the distribution function $\N{\#}$ equals $\Si{\#}$. 
\end{theorem}
For site percolation, the theorem is proven in \cite{Veselic-05a}, 
the modification to bond percolation is not hard. 
See also \cite{KirschM-06} for the case of bond-percolation on the  lattice.

Thus we have obtained an selfaveraging IDS. The theorem implies that zero is the lower edge of the support of $\N{\#}$ for ${\scriptstyle\#}\in
\left\{A,D,N\right\}$ and that $\N{A}(0)=\N{D}(0)=0$ while $\N{N}(0)>0$. 
The distribution function $\N{\#}$ has a discontinuity at $E \in \mathbb{R}$ 
if and only if $E\in \Si{\#}_{fin}$. This will be discussed in more detail in \S \ref{ss-jumps}.
\smallskip

Let us denote the density of vertices in $V \setminus V(\omega)$ by $g(\omega)$.
By ergodicity there is a value $g\in [0,1]$ such that $g(\omega) = g$ almost surely.
In the site percolation model $g= 1-p$ and in the bond percolation model $g=(1-p)^{k}$.
Then the IDS $\widetilde{\N{\#}}$ of $(\widetilde{\D{\#}_\omega})_\omega$
is related to $\N{\#}$ by
 \[
 \widetilde{\N{\#}}= \N{\#} +g\, \chi_{[K,\infty[}.
 \]

\begin{remark}\label{rem: operator comparison}
Let $\mu$ be a non-trivial probability measure on a compact interval $[0,s]$ and
$(w_\omega(x))_{x\in V}$  an i.i.d.~family of random variables with distribution $\mu$.
The Anderson model on the graph $G$ is given by $H_\omega:=\Delta + W_\omega\colon \ell^2(V)\to\ell^2(V)$.
Here $W_\omega$ is the multiplication operator with $w_\omega(x)$. This is again an ergodic operator with 
non-random spectrum and selfaveraging IDS. It is very convenient to compare it with the site percolation Laplacian.
For this purpose we specialise to the case $\mu=p \delta_0 + (1-p) \delta_1$ and introduce an additional, 
global coupling constant $\lambda\ge 0$. We denote the quadratic form of the 
operator $H_{\omega,\lambda}^{\scriptscriptstyle BA}:=\Delta + \lambda W_\omega$ by $\q{BA}_{\omega,\lambda}$. 
Here $BA$ stands for the \emph{Bernoulli-Anderson} or \emph{binary-alloy} type of the Hamiltonian.
Denote its IDS by $\N{BA}_\lambda$.

Recall that one may compare two (lower-bounded, closed) quadratic forms $q_1$ and $q_2$ on a Hilbert space
by setting $q_1 \le q_2$ if and only if the domains satisfy the inclusion $\cD(q_1)\supset \cD(q_2)$ and
for all $\varphi\in \cD(q_2)$ we have $q_1[\varphi] \le q_2[\varphi] $.
This notion of inequality of quadratic forms is consistent with extending each form to the whole
Hilbert space by the rule $q_i[\varphi]=+\infty$ for $\varphi \not\in \cD(q_i)$.

Denote the quadratic form associated to any of the Hamiltonians $\D{\#}_\omega$ 
by $\q{\#}_\omega$.
Comparing the Anderson and the adjacency site percolation Hamiltonian we see that 
$\q{BA}_{\omega,\lambda} \le \q{A}_\omega$ for all $\lambda\ge 0$. This relation holds pointwise 
for all $\omega$ if we introduce the obvious coupling between the random potential
$(w_\omega(x))_{x\in V}$ and the percolation process. Moreover for $\lambda \to \infty$
the quadratic form $\q{BA}_{\omega,\lambda}$ converges monotonously to $\q{A}_\omega$.
 On the other hand, since the potential $(w_\omega(x))_{x\in V}$  is non-negative
we have $q_G\le \q{BA}_{\omega, \lambda}$ where $q_G$ denotes the quadratic form of 
$\Delta$ on the full graph. 
Using the argument in \cite[\S 2]{Simon-85b} one sees that $\widetilde{\q{N}_\omega} \leq q_G$, 
where $\widetilde{\q{N}_\omega}$ is the form corresponding to the operator 
$\widetilde{\D{N}_\omega}= \D{N}_\omega\oplus 0 \cdot\Id_{V \setminus V(\omega)} $.

Summarising all quadratic form inequalities obtained so far for the site percolation case
we obtain the chain 
 \[
 \widetilde{\q{N}_\omega} \leq q_G \le \q{BA}_{\omega,\lambda} \le \q{A}_\omega \le \q{D}_\omega
 \]
which implies for the corresponding IDS'
\[ 
 \widetilde{\N{N}}(E) \ge N_G(E) \ge \N{BA}_\lambda(E) \ge \N{A}(E) \ge \N{D}(E)
\ \text{ for all } E\in\RR.
\]
Now it is clear why the study of the Anderson Hamiltonian in \cite{BiskupK-01a} 
gave also results on the (adjacency) percolation Laplacian.
Upper bounds on the IDS $\N{BA}_\lambda$
imply upper estimates for the IDS of the adjacency and Dirichlet site percolation model.

There is also a relation between bond percolation Hamiltonians on the one hand and random hopping models 
\cite{KloppN-03} and discrete Schr\"odinger operators with random magnetic field \cite{Nakamura-00} on the other hand.
All three models share the feature that the randomness enters the Hamiltonian in the off-diagonal matrix elements.
In fact, in \cite{Nakamura-00,KloppN-03} certain spectral properties of such models are analysed by 
relating them to Hamiltonians with diagonal disorder.
\end{remark}

\begin{remark} \label{rem: bands}
The Anderson model on $\ell^2(\ZZ^d)$ is interpreted by physicists
as a single spectral band approximation of a continuum random Schr\"odinger operator on $L^2(\RR^d)$.
Due to the i.i.d.~assumption on the random variables in the Anderson model 
its spectrum as a set equals almost surely
\[
\sigma(\Delta) + \supp \mu := \{t + s \in \RR \mid t \in \sigma(\Delta), s \in \supp \mu  \}
\]
see Theorem III.2 in \cite{KunzS-80}.
Here as before $\mu $ denotes the distribution measure of the potential values of the 
Anderson model and $\supp \mu$ its (topological) support. Now, if the support of $\mu$ is split into several connected components
and if the gaps between them are large enough, the a.s.~spectrum  of the Anderson model
also contains gaps and consequently has a (sub)band  structure. 
Internal Lifshitz tails for such models have been studied in \cite{Mezincescu-86,Simon-87}.

In particular, if we consider the Bernoulli-Anderson model $H_{\omega,\lambda}^{\scriptscriptstyle BA}$ introduced above, 
we see that for large enough $\lambda$ the almost sure spectrum of  $H_{\omega,\lambda}^{\scriptscriptstyle BA}$
splits into subbands. As we let $\lambda\to\infty$, one of the bands diverges to infinity. 
On the other hand, we know that in the sense of quadratic forms 
$\q{BA}_{\omega,\lambda} \nearrow \q{A}_\omega$  for $\lambda\to\infty$.
Thus the resulting site percolation Hamiltonian $\D{A}_\omega$ may be understood as an approximation of  
$H_{\omega,\lambda}^{\scriptscriptstyle BA}$ (with large values of $\lambda$)  associated to a spectral (sub)band.
\end{remark}

Now we are in the position to state our main results in the next two theorems.
In both cases we have the same setting as before: 
we consider a Cayley graph $G$ of an amenable finitely generated group $\Gamma$, 
site or bond percolation on $G$ and the percolation Laplacians $\D{\#}_\omega$ with associated IDS $\N{\#}$.
We restrict ourselves now to the subcritical phase $p < p_c$.
The asymptotic behaviour of the IDS of the adjacency and the Dirichlet percolation Laplacian at low energies is as follows:

\begin{theorem}\label{thm: main Dirichlet}
Assume that $G$ has polynomial growth and $V(n) \sim n^{d}$. Then there are positive constants
$\alpha_{D}^{+}(p)$ and $\alpha_{D}^{-}(p)$ such that for all positive $E$ small enough
\begin{equation}\label{eq: main dirichlet poly}
e^{-\alpha_{D}^{-}(p)E^{-d/2}} \leq \N{D}(E) \leq \N{A}(E) \leq e^{-\alpha_{D}^{+}(p) E^{-d/2}}.
\end{equation}
Assume that $G$ has superpolynomial growth. Then
\begin{equation}\label{eq: main Dirichlet super}
\lim_{E \searrow 0} \frac{\ln |\ln \, \N{D}(E)|}{|\ln E|} 
= \lim_{E \searrow 0} \frac{\ln |\ln \, \N{A}(E)|}{|\ln E|} = \infty.
\end{equation}
\end{theorem}

In particular, the IDS is very sparse near $E=0$ and consequently the bottom of the spectrum is a fluctuation boundary.
Relations \eqref{e-vanHove} and \eqref{eq: main dirichlet poly} imply that in the case of polynomial growth the 
Lifshitz exponent coincides with the van Hove exponent of the Laplacian on the full Cayley graph.
In the case of superpolynomial growth we have that both exponents are infinite. 
One may ask whether the limits \eqref{e-vanHove infinite} and \eqref{eq: main Dirichlet super}, 
which define the two exponents, diverge at the same rate. To express this in a quantitative way 
let us note that in the case of polynomial growth we have 
\begin{equation*}
\lim_{E \searrow 0} \frac{\ln|\ln \, \N{A}(E)| }{|\ln N_0(E)|}\   =1 
\end{equation*}
and the analogous relation for the Dirichlet Laplacian. 
In the case of superpolynomial growth one may hope to prove in certain cases
\begin{equation*}
\lim_{E \searrow 0} \frac{\ln\ln|\ln \, \N{A}(E)| }{\ln|\ln N_0(E)|}\   =1 .
\end{equation*}
We are not able prove that in general, 
but at least for the case of the Lamplighter group, see \S \ref{ss-Lamplighter}.

In the case of Neumann boundary conditions we know that $\widetilde{\D{N}_\omega} \leq \Delta $.
Consequently, the IDS $\widetilde{\N{N}}$ is at the bottom of the spectrum at least as `thick' as the one of the free Laplacian,
which exhibits a van Hove singularity. In particular, the energy zero is not a fluctuation boundary.
This remains true if we pass from $\widetilde{\N{N}}$ to $\N{N}$ by removing the point mass $g\delta_0$ 
from the density of states  measure.

\begin{theorem}\label{thm: main Neumann}
There exist positive constants $\alpha_{N}^{+}(p)$ and $\alpha_{N}^{-}(p)$ such that for all positive $E$
small enough
\begin{equation}\label{eq: main neumann}
e^{-\alpha_{N}^{-}(p)E^{-1/2}} \leq \N{N}(E)-\N{N}(0) \leq e^{-\alpha_{N}^{+}(p) E^{-1/2}}.
\end{equation}
\end{theorem}

The value $\N{N}(0)$ coincides with the average number of clusters per vertex in the random graph $G(\omega)$.
After subtracting this value we can speak of \eqref{eq: main neumann} as a kind of `renormalised'
Lifshitz asymptotics with exponent $1/2$.

\begin{remark}
(a) 
In the next section we give a sketch of the proof of the theorems, while the full version will appear elsewhere.
\\
(b)
The lower bounds in \eqref{eq: main dirichlet poly} and \eqref{eq: main neumann} are actually true for all values of $p$.
\\
(c)
In the special case of bipartite Cayley graphs there is a relation between 
the behaviour of the IDS near the upper and lower spectral edges of the spectrum, see Remark 
\ref{rem: properties}. In that case we are able to characterise the asymptotic behaviour 
of the IDS at the upper spectral boundary in the same way as done in \cite{KirschM-06}.

\end{remark}

\section{Discussion, additional results and sketch of proofs}

In this final section we conclude the paper with a discussion of 
further spectral properties of percolation Hamiltonians, an outline of the proof of 
the Theorems \ref{thm: main Dirichlet} and \ref{thm: main Neumann}, and an important example. 
The example concerns the Lamplighter group, which is amenable, but of exponential growth.
The spectral properties we mentioned are related to jumps of the IDS, to finitely supported eigenfunctions 
and to the \emph{unique continuation principle}, see  \cite{Veselic-05b,Veselic-06} and \cite{LenzV}.

\subsection{Discontinuities of the IDS\label{ss-jumps}}

Let us recall a result on the location of the set of discontinuities of the IDS,
established in \cite{Veselic-05b} for site percolation Laplacians and generalisations thereof.
The result and its proof apply verbatim for the bond percolation model. Furthermore for this 
result it is not necessary to assume that the percolation process is independent at different sites, 
but merely that it is ergodic with respect to the transformations $(\tau_\gamma)_{\gamma  \in \Gamma}$ 
introduced in \S \ref{ss-classical-percolation}.

\begin{proposition}
\label{p-equiv}
Let $G$ be a Cayley graph of an amenable group $\Gamma$
and $(\D{\#}_\omega)_\omega$, ${\scriptstyle \#} \in\{A,D,N\}$
a percolation Laplacian associated either to site or to bond percolation. 
Then the  following two properties are equivalent:
\begin{enumerate}[\rm(i)]
\item
the IDS of $(\D{\#}_\omega)_\omega$  is discontinuous at $E\in\RR$,
\item
$E\in \Si{\#}_{fin}$.
\end{enumerate}
\end{proposition}

\begin{remark}
This result can be strengthened to describe the size of the jumps also.
In a forthcoming paper \cite{LenzV} of Lenz and the second named author
ergodic Hamiltonians on discrete structures are analysed. 
The abstract setting considered there covers in particular Anderson Hamiltonians and
site and bond percolation Laplacians $(\D{\#}_\omega)_\omega$  
on a Cayley graph $G$ of an amenable group $\Gamma$.
Again one has to assume that the stochastic process which determines 
the random potential, respectively the percolation process entering the 
percolation Laplacian, is ergodic with respect to the transformation group
$(\tau_\gamma)_{\gamma  \in \Gamma}$. Then the eigenspace corresponding 
to the eigenvalue $E$ is spanned by finitely supported eigenfunctions
almost surely. See also \cite{Kuchment-91,Kuchment-05} for a similar result for periodic operators
on graphs with an abelian group structure.
\end{remark}


The jumps of the IDS of invariant operators on Cayley graphs 
have been studied in the literature on $L^2$-invariants. The $k$-th $L^2$-Betti number
is the trace per unit volume of the kernel of the Laplacian on $k$-forms. 
It has been proven in great generality that it does not matter whether one defines it via the 
continuum or the combinatorial Laplacian, see e.g.~\cite{Dodziuk-77}. Likewise it is known that
the zeroth $L^2$-Betti number vanishes on amenable Cayley graphs \cite[Thm.~1.7]{Lueck-01},
which is the same as saying that the IDS is continuous at the bottom of the spectrum.
Using the same terminology for random  ergodic operators, we may say that the 
Neumann percolation Laplacian has non vanishing (zero order) $L^2$-Betti number.
With this regard it would be interesting to find an interpretation of higher 
order Betti numbers in terms of quantities of mathematical physics.
The sizes of jumps of the IDS at discontinuity points have be discussed 
for Abelian periodic models in \cite{Sunada-88,Donnelly-81}. The task of characterising 
the sizes of jumps of the IDS is related to the Atiyah conjecture, 
see e.g.~\cite{GrigorchukLSZ-00,Lueck-02} or \cite[Conjecture 2.1]{Lueck-01},
and the references therein.

\subsection{Outline of the proof of Theorems \ref{thm: main Dirichlet} and \ref{thm: main Neumann}}
While the full proofs of our results will be given in \cite{AntunovicV-c} we
present here certain key estimates.
In particular, we state upper and lower  bounds for the IDS of percolation Laplacians in a neighbourhood of the lower
spectral edge. The bounds are given in terms of $\l{\#}(G')$, the lowest non-zero eigenvalue of $\D{\#}(G')$.
These estimates are generalisations of Lemmata 2.7 and 2.9 in \cite{KirschM-06}. The proof of
Theorem \ref{thm: upper bounds for IDS} uses the exponential decay from Theorem \ref{thm: exponential decay}.

\begin{theorem} \label{thm: upper bounds for IDS}
Let $G$ be an amenable Cayley graph and ${\scriptstyle\#}\in \left\{A,D,N\right\}$. 
Assume that there is a continuous strictly decreasing function $f \colon \left[1,\infty\right[ \to \mathbb{R}^{+}$ 
such that $\lim_{s \to\infty}f(s)=0$ and $\l{\#}(G') \geq f(|G'|)$ for any finite subgraph $G'$. 
Then, for every $0<p<p_{c}$ the IDS satisfies the following inequality
\begin{equation} \label{eq: upper bounds for IDS 1}
\N{\#}(E)-\N{\#}(0) \leq e^{-a_{p}f^{-1}(E)},
\end{equation}
for some positive constant $a_{p}$, when $E$ is small enough. Here  $f^{-1}$ denotes the inverse function of $f$.
\end{theorem}

\begin{theorem} \label{thm: lower bounds for IDS}
Let $G$ be an amenable Cayley graph and ${\scriptstyle\#}\in \left\{A,D,N\right\}$. Suppose that there is a
sequence of connected subgraphs $(G'_{n})_{n}$ and a sequence $(c_{n})_{n}$ in $\mathbb{R}^{+}$ such that
\begin{enumerate}[\rm(i)]
    \item $\lim_{n \to \infty} |G'_{n}| = \infty$
    \item $\lim_{n \to \infty}c_{n} = 0$
    \item $\l{\#}(G'_{n}) \leq c_{n}$
\end{enumerate}
For every $E>0$ small enough define $n(E):=\min\left\{n; c_{n} \leq E \right\}$. Then for every $0<p<1$ there is a
positive constant $b_{p}$ such that the following inequality holds for all $E>0$ small enough
\begin{equation}\label{eq: lower bounds for IDS 1}
\N{\#}(E)-\N{\#}(0) \geq e^{-b_{p}|G'_{n(E)}|}.
\end{equation}
\end{theorem}

Hence our problem is reduced to finding efficient bounds for $\l{\#}(G')$ in terms of the 
geometric properties  of $G'$. 
For this we will use, following \cite{KirschM-06}, the Cheeger and Faber \& Krahn inequalities.
Since we are considering general Cayley graphs of amenable groups with polynomial growth,
these two inequalities are not sufficient, but we will need additionally an appropriate version of the isoperimetric inequality.

The function $V$ was defined in Theorem \ref{thm: growth of groups}. We also define 
$$\phi(t):= \min \left\{n \geq 0;V(n)>t \right\}.$$ 
Moreover we will denote the linear subgraph with $n$ vertices by $L_{n}$.

\begin{proposition}\label{prop: eigenvalue bounds}
For a Cayley graph $G=(V,E)$ there are positive constants $\alpha_{D}$, $\beta_{D}$, $\gamma_{D}$, $\alpha_{N}$ and
$\gamma_{N}$ such that the following are true
\begin{itemize}
\item[(i)] For every finite connected subgraph $G'$
\begin{equation}\label{eq: lower bounds}
\l{A}(G') \geq \frac{\alpha_{D}}{\phi(\beta_{D} |G'|)^{2}} \ \ \textrm{ and } \ \ \l{N}(G')
\geq \frac{\alpha_{N}}{|G'|^{2}}.
\end{equation}
\item[(ii)] For every positive integer $n$
\begin{equation}\label{eq: upper bounds}
\l{D}(B(n)) \leq \frac{\gamma_{D} V(n)}{n^{2} V(\lfloor n/2 \rfloor)} \ \ \textrm{ and } \ \ \l{N}(L_{n})
\leq \frac{\gamma_{N}}{n^{2}}.
\end{equation}
\end{itemize}
\end{proposition}

\begin{proof}[Sketch of the proof] 
\begin{itemize}
\item[(i)] The inequality for adjacency Laplacian can be proved following the arguments in the proofs of Lemma
2.4 from \cite{KirschM-06} and Proposition 7.1 from \cite{ChungGY-00} and
using the isoperimetric inequality from Th\'{e}or\`{e}me 1 in \cite{CoulhonSC-93}.\\
The inequality for Neumann Laplacian is a simple consequence of the Cheeger inequality (see Th\'{e}or\`{e}me 3.1 in
\cite{ColindeVerdiere-98}). \item[(ii)] Both inequalities can be obtained by inserting an appropriate test function to the
sesquilinear form and using the mini-max principle. For the Dirichlet Laplacian we choose a test function which has
value $0$ outside the ball $B(n)$, value $i$ on the sphere of radius $n-i$, for $i=0,\dots,\lceil n/2 \rceil$ and value
$\lceil n/2 \rceil$
inside the ball $B(\lfloor n/2 \rfloor)$.\\
The operator $\D{N}(L_{n})$ is not injective and thus the test function must be orthogonal to the kernel. This
means that the test function $\varphi$ must have support inside $L_{n}$ and $\sum_{x \in L_{n}}\varphi(x) = 0$. The
appropriate test function is the one which grows linearly from $\frac{-n+1}{2}$ to $\frac{n-1}{2}$ along the vertices
of $L_{n}$ and is $0$ outside of $L_{n}$.
\end{itemize}
\end{proof}

\begin{remark}
(a) \,
In the case of polynomial growth the first bound in \eqref{eq: lower bounds} and the first bound in 
\eqref{eq: upper bounds} are of the same order in $n$ if $G'=B(n)$. Likewise the second bounds in 
\eqref{eq: lower bounds} and \eqref{eq: upper bounds} are of the same order of magnitude for $G'=L_n$.
Thus the proposition shows that the optimal subgraph configuration for adjacency and Dirichlet Laplacians  
is a ball, and for the Neumann Laplacian is a line graph. 
For Dirichlet and Neumann boundary conditions this can be motivated in terms of the potential $\pm W^{b.c.}$.
Since $W^{b.c.}$ is non-negative (repulsive) at the interior vertex boundary, a test function $\varphi$ with 
low expectation value $\la \varphi, \D{D}(G') \varphi\ra$  
has to make both $\la \varphi, \D{A}(G')\varphi\ra$ and $\la \varphi, W^{b.c.}\varphi\ra$ 
small. The second condition pushes the mass of $\varphi$  away form the boundary, 
while the first one minimises the variation of $\varphi$. This is best realised when $G'$ is a ball.
On the contrary, in the case of Neumann boundary conditions $-W^{b.c.} \le 0$ is attracting 
$\varphi$ to the boundary. For fixed volume the graph with most boundary is a linear graph.

(b) \,
Note that if $\Gamma$ is a group of polynomial growth, then the balls $B(n), n\in\NN$ form a F\o lner sequence:
Since every finite $F$ satisfies $F\subset B(k)$ for some $k \in \NN$ we have 
\[
|F \cdot B(n)\setminus B(n)| \le V(n+k)-V(n) \le (c+\co(1))  \, (n+k)^d -(c-\co(1)) \, n^d = \co(V(n))
\]
and similarly $|B(n)\setminus F \cdot B(n)| \le V(n)-V(n-k) =  \co(V(n))$. 
Since $V(2n) \le 2^d \frac{b}{a} V(n)$ the sequence of balls satisfies the \emph{doubling property}
and is in particular tempered. Note that although every amenable group contains a F\o lner sequence,
balls may not form one.

\end{remark}

\subsection{The Lamplighter group\label{ss-Lamplighter}}

In this section we will explain how the ideas and methods we used
to study the behaviour of the IDS in the case of groups of
polynomial growth, can be used to give sharp bounds on the IDS
asymptotics in the case of a particular group of superpolynomial
growth. Namely, we consider certain Lamplighter groups, which
are examples of amenable groups with exponential growth
(i.e.~there exists a
constant $c>1$ such that $V(n) \geq c^{n}$).

We define the Lamplighter group as the wreath product
$\mathbb{Z}_{m} \wr \mathbb{Z}$, where $m$ is an arbitrary
positive integer. Elements of this group are ordered pairs
$(\varphi,x)$, where $\varphi$ is a function $\varphi \colon
\mathbb{Z} \to \mathbb{Z}_{m}$ with finite support and $x\in
\mathbb{Z}$. The multiplication is given by $(\varphi_{1},x_{1})
\ast (\varphi_{2},x_{2}) := (\varphi_{1} + \varphi_{2}(\cdot -
x_{1}),x_{1}+x_{2})$.
More generally, a lamplighter group can be defined as a wreath
product of two Abelian groups, one of which is finite. However, here we
shall only consider lamplighter groups of the form $\mathbb{Z}_{m}
\wr \mathbb{Z}$.

Since Lamplighter groups are amenable it makes sense to consider the
integrated density of states in both the deterministic and random
setting, as defined in Section \ref{ss-results}.
The next theorem concerns the Laplace operator on the full Cayley graph. 
It establishes a  Lifshitz-type asymptotics in the sense that 
\begin{equation}
\label{eq: Lamplighter exponent}
\lim_{E\searrow 0} \frac{\ln |\ln N_{0}(E)|}{|\ln E|} = \frac{1}{2}.
\end{equation}
Of course the underlying operator is not random, but the IDS is exponentially thin 
at the bottom of the spectrum as it is the case for spectral fluctuation  boundaries 
of random operators. The quantity \eqref{eq: Lamplighter exponent}
may be called a Lifshitz exponent or maybe more appropriately 
secondary Novikov-Shubin invariant of the Laplacian, using the terminology of \cite{Oguni-07}.

\begin{theorem}\label{thm: lamplighter deterministic}
Let $G$ be a Cayley graph of the Lamplighter group. There are
positive constants $a_{1}^{+}$ and $a_{2}^{+}$ such that
$$
N_{0}(E) \leq a_{1}^{+}e^{-a_{2}^{+}E^{-1/2}}, \textrm{ for all }
E \textrm{ small enough.}
$$
Moreover for every $r>1/2$ there are positive constants
$a_{r,1}^{-}$ and $a_{r,2}^{-}$ such that
$$
N_{0}(E) \geq a_{r,1}^{-}e^{-a_{r,2}^{-}E^{-r}}, \textrm{ for all
} E \textrm{ small enough.}
$$
\end{theorem}
The proof of the preceding theorem follows from the proof of
Theorem 4.4 (parts (ii) and (iii)) in \cite{Oguni-07}. As an input
we need the following inequalities for $\mu^{(2n)}(\iota)$, the
return probability  of the simple random walk after $2n$
steps:
$$
a_{1}e^{-a_{2}n^{1/3}} \leq \mu^{(2n)}(\iota) \leq
A_{1}e^{-A_{2}n^{1/3}},
$$
for some positive constants $a_{i}$, $A_{i}$, $i=1,2$, and all
positive integers $n$. For the reference see Theorem 15.15 in
\cite{Woess-00}.
\smallskip

For the percolation case we shall once again use Theorems
\ref{thm: upper bounds for IDS} and \ref{thm: lower bounds for
IDS}.
Exponential growth and Proposition \ref{prop: eigenvalue bounds}
i) give lower bounds for the lowest eigenvalues $\l{A}(G')$
which are of the form $const/(\ln |G'|)^{2}$. Now Theorem
\ref{thm: upper bounds for IDS} implies the following result.

\begin{theorem}\label{thm: lamplighter upper bounds}
Let $G$ be an arbitrary Cayley graph of the lamplighter group
$\mathbb{Z}_{m} \wr \mathbb{Z}$. For every $p<p_{c}$ there are
positive constants $b_{1}$ and $b_{2}$ such that the IDS of
the adjacency and Dirichlet percolation Laplacian satisfies the following inequality
\begin{equation}\label{eq: lamplighter upper bounds}
N^{A}(E) \leq N^{D}(E) \leq e^{-b_{1}e^{b_{2}E^{-1/2}}}, \textrm{
for all } E \textrm{ small enough.}
\end{equation}
\end{theorem}

The  upper bounds from Proposition \ref{prop: eigenvalue bounds}
are not applicable in the Theorem \ref{thm: lower bounds for IDS}
any more. Instead, we shall use results from \cite{BartholdiW-05}, where
Bartholdi \& Woess proved that the Laplacian $\Delta$ on Diestel-Leader graphs
has only pure point spectrum, expanding on earlier results of Grigorchuk \& {\.Z}uk \cite{GrigorchukZ-01}
and Dicks and Schick \cite{DicksS-02} concerning certain Cayley graphs of Lamplighter groups.
Diestel-Leader graphs are a generalisation of Cayley graphs of Lamplighter groups with 
a particular set of generators. This `natural' set of generators is given by
\begin{equation}\label{eq: generator set}
S_{0}:=\left\{(l\cdot\delta_{1},1), l \in \mathbb{Z}_{m}\right\}
\cup \left\{(l\cdot\delta_{0},-1), l \in \mathbb{Z}_{m}\right\}.
\end{equation}
Here $l\cdot\delta_{z}$ denotes the function which has value $l$ in
$z$ and $0$ everywhere else. 
In the following we formulate certain facts about the spectrum of
adjacency Laplacians on certain finite, connected subgraphs 
of the Lamplighter graph, called \emph{tetrahedrons},
 which were established in \cite{BartholdiW-05}, see also
\cite{GrigorchukZ-01,DicksS-02}.
We will not give a definition of these subgraphs, since it would require a quite comprehensive
description of horocyclic products of homogeneous trees.
Rather, we refer to \cite{BartholdiW-05} for precise definitions 
and background information.

We will need three facts concerning tetrahedrons and
corresponding eigenvalues of adjacency Laplacians:
\begin{enumerate}[(a)]
\item 
the tetrahedron of depth $n$ (denoted by $T_{n}$) has
$(n+1)m^{n}$ vertices, 
\item
$2m(1-\cos \frac{\pi}{n})$ is an
eigenvalue of the operator $\Delta^{A}(T_{n})$, 
\item
there is
an eigenvector of $\Delta^{A}(T_{n})$, corresponding to the
eigenvalue $2m(1-\cos \frac{\pi}{n})$, that has value $0$ on the
inner vertex boundary of $T_{n}$.
\end{enumerate}
Claim a) follows directly from the definition of tetrahedron. For
the proofs of
claims b) and c) see Lemma 2 and Corollary 1 in \cite{BartholdiW-05}.

Now Theorem \ref{thm: lower bounds for IDS}, with $G_{n}'=T_{n}$
and $c_{n}=2m(1-\cos\frac{\pi}{n})$, implies the following bound.

\begin{theorem}\label{thm: lamplighter lower bounds}
Let $G_{S_{0}}$ be the Cayley graph of the lamplighter group
$\mathbb{Z}_{m} \wr \mathbb{Z}$ defined with respect to the set of
generators $S_{0}$. For every $0<p<1$ there are positive constants
$c_{1}$ and $c_{2}$ such that the IDS of percolation Laplacians
satisfies the following inequality
\begin{equation}\label{eq: lamplighter lower bounds}
e^{-c_{1}e^{c_{2}E^{-1/2}}} \leq N^{D}(E) \leq N^{A}(E), \textrm{
for all } E \textrm{ small enough.}
\end{equation}
\end{theorem}
Theorem \ref{thm: lamplighter lower bounds} is in fact valid for
all Cayley graphs of the lamplighter group $\mathbb{Z}_{m} \wr
\mathbb{Z}$ (with constants $c_{1}$ and $c_{2}$ possibly depending
on the choice of the generator set). This can be seen as follows:

Assume we are given a Cayley
graph $G_{S}$ defined with respect to a generator set $S$. A
natural candidate for the sequence of subgraphs $G_{n}'$ in
Theorem \ref{thm: lower bounds for IDS} would be $G_{S}(V_{n})$,
the subgraphs in $G_{S}$ induced by $V_{n}$ ($V_{n}$ being the vertex
set of a tetrahedron with depth $n$). Since this subgraph is not necessarily
connected, to each vertex $x$ in $V_{n}$ we add a ball
(in $G_{S}$) of some large, but fixed radius $R$, centred at $x$.
In this way we get the vertex set $V_{n}^{R}$, which is connected
in $G_{S}$. Because of fact (c) above, the adjacency Laplacian on the
subgraph of $G_{S_{0}}$ induced by the vertex set $V_{n}^{R}$ is
again bounded above by $2m(1-\cos\frac{\pi}{n})$. Now calculations
similar to those in the proof of Theorem 3.2 in \cite{Woess-00}
show that the lowest eigenvalue of the adjacency Laplacian on the
$G_{S}(V_{n}^{R})$ has an upper bound of the form
$const(1-\cos\frac{\pi}{n})$. Thus Theorem \ref{thm: lower bounds
for IDS} with $G_{n}'= G_{S}(V_{n}^{R})$ proves the claim in
Theorem \ref{thm: lamplighter lower bounds} for the Cayley graph
$G_{S}$.

\subsection*{Acknowledgements}
It is a pleasure to thank D.~Lenz, I.~Naki\'c and F.~Sobieczky 
for stimulating discussions related to the topics of this paper.
Visits of the authors to the TU Graz were supported by
the ESF through a Short Visit Grant
and the \"Ostereichischer Fonds zur F\"orderung der wissenschaftlichen Forschung.

\def\cprime{$'$}
\def\polhk#1{\setbox0=\hbox{#1}{\ooalign{\hidewidth
  \lower1.5ex\hbox{`}\hidewidth\crcr\unhbox0}}}

\end{document}